\documentclass[11pt]{article}%
\usepackage{amsmath}
\usepackage{amsfonts}
\usepackage{amssymb}
\usepackage{hyperref}
\usepackage{graphicx}
\usepackage{theorem}
\usepackage{enumerate}
\usepackage{color}
\usepackage{hyperref}%
\providecommand{\U}[1]{\protect\rule{.1in}{.1in}}
\providecommand{\U}[1]{\protect\rule{.1in}{.1in}}

\newtheorem{thm}{Theorem}[section]

\newtheorem{pr}[thm]{Proposition}
\newtheorem{df}[thm]{Definition}
\newtheorem{rmk}[thm]{Remark}

{\theorembodyfont{\upshape}
\newtheorem{examp}[thm]{Example}
}
\numberwithin{equation}{section} 
\hypersetup{
colorlinks=true,
linkcolor=blue,
citecolor=red,
filecolor=magenta,
urlcolor=cyan
}
\begin{document}

\title{Title: }

\begin{center}
NONLINEAR\ AND\ ORBITAL DIRECTIONAL REGULARITIES\ FOR\ SET-VALUED MAPPINGS
AND\ APPLICATIONS TO FIXED POINTS

\bigskip

by

\bigskip

Marius Durea \footnote{{\small Faculty of Mathematics, \textquotedblleft
Alexandru Ioan Cuza\textquotedblright\ University, 700506--Ia\c{s}i, Romania
and \textquotedblleft Octav Mayer\textquotedblright\ Institute of Mathematics,
Ia\c{s}i Branch of Romanian Academy, 700505--Ia\c{s}i, Romania; e-mail:
\texttt{durea@uaic.ro}}}, Nguyen Huu Tron\footnote{{\small Department of
Mathematics and Statistics, Quy Nhon University, 170 An Duong Vuong, Quy Nhon,
Gia Lai, Vietnam; e-mail: \texttt{nguyenhuutron@qnu.edu.vn}}}, Michel
Th\'{e}ra\footnote{{\small XLIM UMR-CNRS 7252, Universit\'{e} de Limoges,
Limoges, France, ORCID 0000-0001-9022-6406; e-mail:
\texttt{michel.thera@unilim.fr}}}
\end{center}

\bigskip

\noindent Abstract: Motivated by developments in metric regularity theory for
set-valued mappings, this paper introduces several directional regularity
notions that are weaker than the corresponding classical regularity
properties. We investigate the relationships among these concepts and
establish a classification framework based on directional minimal time
functions. Particular attention is devoted to directional orbital regularity
and directional orbital Aubin continuity, which are shown to provide a natural
setting for the analysis of fixed point phenomena. 

Using these notions, we derive approximate and exact fixed point theorems for
set-valued mappings under directional assumptions. The obtained results are
expressed in terms of minimal time functions and orbital approximation
procedures, allowing for anisotropic and asymmetric behaviors that are not
captured by the classical nondirectional framework. We further apply the
developed theory to directional coincidence point results, coupled fixed point
theorems, and Milyutin-type perturbation stability properties. 

The proposed approach provides a unified directional perspective on fixed
point theory and variational analysis. In particular, it recovers several
known results from the literature as special cases and yields new extensions
under substantially weaker assumptions than classical contraction or global
regularity conditions. 

\bigskip



\section{Introduction}

Fixed point theory plays a fundamental role in nonlinear analysis and has
numerous applications in optimization, variational analysis, equilibrium
theory, economics, engineering, data science, and machine learning; see
\cite{CP} and the references therein. Since the pioneering Banach contraction
principle \cite{B,C,GD,KK}, fixed point methods have become indispensable
tools for the analysis of nonlinear equations, differential and integral
equations, variational inequalities, and optimization problems. 

The development of modern applications has led to the study of fixed point
phenomena in increasingly general settings involving set-valued mappings and
generalized metric structures. A cornerstone in this direction is Nadler's
theorem \cite{N} for contractive multifunctions. At the same time, deep links
between fixed point theory and metric regularity have emerged through the
works of Dontchev and Rockafellar \cite[Theorem~5E.2]{DR-Book} and Ioffe
\cite{I1,I2,I3}. These connections have proved particularly fruitful in
variational analysis and generalized equations. 

A major trend in recent years has been the search for weaker regularity
assumptions capable of preserving fixed point and stability properties. In
this context, Ioffe introduced orbital regularity \cite{I3}, a weaker form of
metric regularity adapted to iterative constructions. Later, Ait Mansour,
Bahraoui, and El Bekkali \cite{ABB} established approximate fixed point
results in incomplete metric spaces, and the approaches of \cite{ABB} and
\cite{I3} were unified in \cite{T}. Related developments concerning
coincidence points, cyclic fixed points, and stability theory may be found in
\cite{A1,AAZ,DF,DS,DR-Book,I1,I2,I3,I4,L,NTH2,TT1,TT2,TT3}. 

Parallel to these developments, directional methods have become an active area
of research in variational analysis and nonsmooth optimization. Directional
regularity notions provide a refined description of multifunction behavior by
focusing on prescribed directions rather than all perturbations. Such concepts
are particularly useful when classical metric regularity fails or when the
underlying phenomena exhibit anisotropic features. Various forms of
directional regularity have been investigated in
\cite{BDSTRU,CDS,CMNW,DPS1,DPS2,Gfr1,Gfr2,ZO}. 

A key tool in this theory is the directional minimal time function introduced
and studied in \cite{DPS1,NZ}. Unlike the classical distance function, minimal
time functions naturally encode directional information and have proved
effective in vector optimization, location theory, variational systems, and
generalized regularity analysis. In particular, several directional regularity
properties can be characterized through suitable minimal time functions; see
\cite{CDS,DPS2}. 

The purpose of this paper is to develop a directional regularity framework for
the study of fixed point phenomena associated with set-valued mappings. We
introduce directional versions of Milyutin regularity and orbital regularity
on sets and establish for each notion a complete classification theorem in
terms of directional openness, directional Aubin continuity, and directional
regularity properties. The proposed constructions are generated by subsets of
the unit spheres of the underlying spaces and are naturally formulated through
directional minimal time functions. Classical regularity notions are recovered
as particular cases when the admissible direction sets coincide with the whole
unit spheres. 

Our main objective is to investigate fixed point phenomena generated by
directional orbital regularity and directional orbital Aubin continuity.
Combining directional minimal time functions with orbital approximation
techniques, we establish approximate and exact fixed point theorems under weak
directional assumptions. A central role is played by directional orbital
trajectories generated by successive approximation procedures. The obtained
conditions are substantially weaker than classical contraction assumptions and
naturally accommodate directional asymmetric behaviors. 

The developed theory is subsequently applied in two directions. First, we
derive directional coincidence point and coupled fixed point theorems by means
of suitable product-space constructions. Second, we obtain approximate and
exact Milyutin-type perturbation stability results as consequences of the
directional regularity framework. Our results recover, unify, and
substantially extend several recent contributions from the literature,
including those in \cite{ABB,I3,T}. 

The paper is organized as follows. Section~2 contains the necessary background
material and introduces the directional regularity notions studied throughout
the paper. Section~3 establishes approximate and exact fixed point theorems
via directional orbital regularity in both complete and incomplete metric
spaces. The final main section is devoted to applications, including
directional coincidence point theorems, coupled fixed point theorems, and
directional versions of Milyutin-type stability results. The paper concludes
with several remarks and perspectives for future research. 

\section{Definitions and Preliminary Results}

Let $(X,\Vert\cdot\Vert)$ be a normed vector space. The unit sphere of $X$ is
denoted by $S_{\mathbb{X}}$. For $x\in X$ and $r>0$, we denote by $B(x,r)$ the
open ball of radius $r$ centered at $x$. 

For a nonempty set $\Omega\subset X$, the symbols $\operatorname{cl}\Omega$,
$\operatorname{int}\Omega$, and $\operatorname{cone}\Omega$ stand for the
closure, interior, and conic hull of $\Omega$, respectively. The distance from
$x\in X$ to $\Omega$ is defined by
\[
\mathbf{d}(x,\Omega):=\inf_{u\in\Omega}\Vert x-u\Vert.
\]
Given $\varepsilon\geq0$, the $\varepsilon$-enlargement of $\Omega$ is
\[
\Omega(\varepsilon):=\{x\in\mathbb{X}:\mathbf{d}(x,\Omega)\leq\varepsilon\}.
\]

Let $F:X\rightrightarrows X$ be a set-valued mapping. A point $x\in X$ is
called an $\varepsilon$-fixed point of $F$ if
\[
\mathbf{d}(x,F(x))\leq\varepsilon,
\]
where, as usual, $d(x,\emptyset):=+\infty$. If $\varepsilon=0$ and $F$ is
closed-valued, this reduces to the classical fixed point condition $x\in
F(x)$. 

The sets of $\varepsilon$-fixed points and fixed points of $F$ are denoted by
\[
\operatorname{Fix}(F)(\varepsilon):=\{x\in\mathbb{X}:\mathbf{d}(x,F(x))\leq
\varepsilon\},
\]
and
\[
\operatorname{Fix}(F):=\{x\in\mathbb{X}:x\in F(x)\},
\]
respectively. 

The following notion will play a central role throughout the paper. 

\begin{df}
Let $F:X\rightrightarrows X$ and let $\emptyset\neq L\subset S_{\mathbb{X}}$.
A sequence $(x_{k})_{k\in\mathbb{N}}\subset X$ is called a sequence of
successive approximations of $F$ starting from $\bar{x}\in X$ with respect to
(w.r.t., for short) $L$ if there exist sequences $(t_{k})\subset R$ and
$(u_{k})\subset L$ such that
\[
x_{0}=\bar{x},\qquad t_{k}\rightarrow0,
\]
and
\[
x_{k+1}=x_{k}+t_{k}u_{k}\in F(x_{k})\qquad\forall k\in\mathbb{N}.
\]

\end{df}

A further key ingredient of our analysis is the directional minimal time
function introduced in \cite{NZ} and subsequently developed in
\cite{DPS1,DPS2,CDS}. 

\begin{df}
Let $A\subset X$ and $\emptyset\neq L\subset S_{\mathbb{X}}$. The function
\[
T_{L}(x,A):=\inf\{t\geq0\mid(x+tL)\cap A\neq\emptyset\}
\]
is called the directional minimal time function associated with $L$. By
convention,
\[
T_{L}(x,\emptyset):=+\infty\qquad\forall x\in\mathbb{X}.
\]

\end{df}

It follows directly from the definition that
\[
T_{L}(x,A)<+\infty\iff x\in A-\operatorname{cone}L,
\]
and
\[
\mathbf{d}(x,A)\leq T_{L}(x,A)\qquad\forall x\in\mathbb{X}.
\]
Moreover, if $L=S_{\mathbb{X}}$, then
\[
T_{L}(\cdot,A)=\mathbf{d}(\cdot,A).
\]

For singleton sets $A=\{u\}$, we simply write $T_{L}(\cdot,u)$. In this case,
\[
T_{L}(x,u)<+\infty\iff u-x\in\operatorname{cone}L,
\]
and, whenever this condition holds,
\[
T_{L}(x,u)=\Vert x-u\Vert.
\]
Furthermore,
\[
T_{L}(x,u)=T_{-L}(u,x)\qquad\forall x,u\in\mathbb{X}.
\]
Moreover, one also has (see \cite[Lemma 1.2]{BDSTRU}) that
\[
T_{L}(x,\Omega)=\inf_{w\in\Omega}T_{L}(x,w).
\]

If $\operatorname{cone}L$ is convex, then
\[
T_{L}(x,u)=0\iff x=u,
\]
and the triangle inequality
\[
T_{L}(x,u)\leq T_{L}(x,v)+T_{L}(v,u)
\]
holds for all $x,v,u\in X$. 

The directional excess associated with $L$ is defined by
\[
\mathbf{e}_{L}(A,B):=\sup_{x\in A}T_{L}(x,B).
\]
Clearly,
\[
\mathbf{e}_{L}(A,B)=+\infty\qquad\text{whenever }A\not \subset
B-\operatorname{cone}L.
\]
In the particular case $L=S_{\mathbb{X}}$, one recovers the classical excess
\[
e(A,B):=\sup_{x\in A}\mathbf{d}(x,B).
\]

Throughout the paper, we adopt the conventions
\[
\mathbf{e}_{L}(\emptyset,B):=0\quad(B\neq\emptyset),
\]%
\[
\mathbf{e}_{L}(A,\emptyset):=+\infty,
\]
and
\[
\mathbf{e}_{L}(\emptyset,\emptyset):=+\infty.
\]
Let $(Y,\Vert\cdot\Vert)$ be a normed space. For nonempty sets $\Omega
_{1},\Omega_{2}\subset X$ and $\Omega_{3}\subset Y$, and for $L\subset
S_{\mathbb{X}}$, $M\subset S_{\mathbb{Y}}$, one has for all $(x,y)\in X\times
Y$,
\[
T_{L\times M}\bigl((x,y),\Omega_{1}\times\Omega_{3}\bigr)=\max\{T_{L}%
(x,\Omega_{1}),\,T_{M}(y,\Omega_{3})\},
\]
and
\[
T_{L}(x,\Omega_{1}\cup\Omega_{2})=\min\{T_{L}(x,\Omega_{1}),\,T_{L}%
(x,\Omega_{2})\}.
\]

Directional metric regularity notions for set-valued mappings around points of
their graphs were introduced and studied in \cite{DPS2,CDS}. In the present
work, we develop nonlocal and nonlinear extensions based on regularity with
respect to fixed sets, in the spirit of recent developments in
\cite{I1,I2,I3,T,THN}. We also note alternative approaches to metric
(sub)regularity in \cite{Gfr1,Gfr2}. 

\medskip

Throughout, for a set-valued mapping $F:X\rightrightarrows Y$,
$\operatorname{Gr}F$ denotes its graph. A function $\gamma:X\rightarrow
R_{+}\cup\{+\infty\}$ is called a gauge function on $U\subset X$ if
$\gamma(x)>0$ for all $x\in U$. A canonical example is $\gamma
(x)=d(x,X\setminus U)$. 

\begin{df}
\label{Nonlocal} Let $F:X\rightrightarrows Y$, let $\emptyset\neq L\subset
S_{\mathbb{X}}$, $\emptyset\neq M\subset S_{\mathbb{Y}}$, $U\subset X$,
$V\subset Y$, and let $\gamma,\delta$ be gauge functions on $U$ and $V$,
respectively. Let $\tau>0$. The mapping $F$ is 

\begin{enumerate}
\item[(i)] directionally $\gamma$-regular on $U\times V$ w.r.t.\ $L\ and$ $M$
with constant $\tau$ if for every $(x,y)\in U\times V$
\[
\tau T_{M}(y,F(x))<\gamma(x)\ \Longrightarrow\ T_{L}(x,F^{-1}(y))\leq\tau
T_{M}(y,F(x)).
\]

\item[(ii)] directionally $\delta$-Aubin continuous on $U\times V$
w.r.t.\ $L\ and$ $M$ with constant $\tau$ if for every $x\in U,$ $u\in X$ and
$y\in F\left(  u\right)  \cap V$
\[
\tau T_{L}(x,u)<\delta(y)\ \Longrightarrow\ T_{M}(y,F(x))\leq\tau T_{L}(x,u).
\]

\item[(iii)] directionally $\gamma$-open on $U\times V$ w.r.t.\ $L\ and$ $M$
with constant $\tau$ if for every $(x,y)\in(U\times Y)\cap\operatorname{Gr}F$
and every $t\in(0,\gamma(x))$,
\[
\mathbf{B}(y,\tau^{-1}t)\cap V\cap\lbrack y-\operatorname{cone}M]\subset
F\bigl(\mathbf{B}(x,t)\cap\lbrack x+\operatorname{cone}L]\bigr).
\]

\end{enumerate}
\end{df}

\begin{rmk}
\begin{enumerate}

\item[(i)] If $L=S_{\mathbb{X}}$ and $M=S_{\mathbb{Y}}$, the above notions
reduce to the classical nonlinear regularity concepts for sets; see
\cite{I1,I3}. Constant gauge functions recover the standard (local) regularity
framework (cf. \cite{DS}). 

\item[(ii) ] The case $U=X$ and $V=Y$ corresponds to global properties. 

\item[(iii)] Choosing $U=B(\bar{x},\varepsilon)$, $V=B(\bar{y},\varepsilon)$
with $(\bar{x},\bar{y})\in\operatorname{Gr}F$ and constant $\gamma
\equiv\varepsilon$ recovers the local framework of \cite{DPS2}. 

\item[(iv)] Taking $\gamma(x)=d(x,X\setminus U)$ and $\delta(y)=d(y,Y\setminus
V)$ yields the Milyutin-type regularity notions; see \cite[Def.~2.28]{I3}. 
\end{enumerate}
\end{rmk}

The next result is from \cite{ZO}. 

\begin{examp}
Let $X=Y=R$, $L=M=\{1\}$, and define
\[
F(x)=%
\begin{cases}
\{x\}, & x\geq0,\\
\emptyset, & x<0,
\end{cases}
\qquad G(x)=%
\begin{cases}
\{\tfrac{1}{2}x\}, & x\geq0,\\
\emptyset, & x<0.
\end{cases}
\]
With constant gauges $\gamma,\delta$, the mapping $F$ is directionally
$\gamma$-regular at $(0,0)$ (w.r.t. $L$ and $-M)$ and constant $1$), while $G$
is directionally $\delta$-Aubin continuous at $(0,0)$ with modulus $1/2$. 

However, $F$ fails to be metrically regular at $(0,0)$, and $G$ is not Aubin
continuous in the classical sense (corresponding to $L=M=\{-1,1\}$). 
\end{examp}

Further examples can be found in \cite[Ex.~2.5]{DPS1} and \cite[Ex.~6]{CDS}. 

\bigskip

The above directional regularity, Aubin continuity, and openness properties
are closely intertwined, mirroring the classical triad of regularity notions;
see \cite{THN,I3,DPS2}. Their interrelations are established in the next
result. 

\begin{pr}
\label{Equivalence} Let $F:X\rightrightarrows Y$ be a set-valued mapping. Let
$L\subset S_{\mathbb{X}}$, $M\subset S_{\mathbb{Y}}$, and let $U\subset X$,
$V\subset Y$ be nonempty sets. Let $\gamma$ be a gauge function on $U$, and
let $\tau>0$. Then the following assertions are equivalent: 

\begin{enumerate}
\item[(i)] $F$ is directionally $\gamma$-regular on $U\times V$ w.r.t. $L$ and
$M$ with constant $\tau$; 

\item[(ii)] $F^{-1}$ is directionally $\gamma$-Aubin continuous on $V\times U$
w.r.t. $M$ and $L$ with constant $\tau$; 

\item[(iii)] $F$ is directionally $\gamma$-open on $U\times V$ w.r.t. $L$ and
$M$ with constant $\tau^{-1}$. 
\end{enumerate}
\end{pr}

\noindent Proof.  We start by proving that $(i)\implies(ii).$ Consequently, we
take $y\in V,v\in Y$ and $x\in F^{-1}\left(  v\right)  \cap U$ with $\tau
T_{M}\left(  y,v\right)  <\gamma\left(  x\right)  .$ Since $v\in F\left(
x\right)  $ one has that $\tau T_{M}\left(  y,F\left(  x\right)  \right)
<\gamma\left(  x\right)  $ and by our assumption $T_{L}\left(  x,F^{-1}\left(
y\right)  \right)  \leq\tau T_{M}\left(  y,F\left(  x\right)  \right)
\leq\tau T_{M}\left(  y,v\right)  .$ 

We next prove $(ii)\implies(iii)$. Let $(x,y)\in\operatorname{Gr}F\cap(U\times
Y)$ and let $t\in(0,\gamma(x))$. Consider
\[
v\in\mathbf{B}(y,\tau^{-1}t)\cap V\cap\lbrack y-\operatorname{cone}M].
\]
Then
\[
\tau T_{M}(v,y)<t.
\]
Since $x\in F^{-1}(y)$, assumption (ii) yields
\[
T_{L}(x,F^{-1}(v))<t.
\]
Hence there exists
\[
u\in F^{-1}(v)\cap\mathbf{B}(x,t)\cap\lbrack x+\operatorname{cone}L].
\]
Therefore,
\[
v\in F\bigl(\mathbf{B}(x,t)\cap\lbrack x+\operatorname{cone}L]\bigr),
\]
which proves $(iii)$. 

Finally, we prove that $\left(  iii\right)  \implies(i).$ Consider $(x,y)\in
U\times V$ satisfying $\tau T_{M}\left(  y,F(x)\right)  <\gamma(x).$ Then for
all $\varepsilon>0$ one can find $v_{\varepsilon}\in F\left(  x\right)
\cap\left[  y+\operatorname*{cone}M\right]  $ such that $\tau\left\Vert
v_{\varepsilon}-y\right\Vert <\tau\left(  T_{M}\left(  y,F(x)\right)
+\varepsilon\right)  .$ Of course, for $\varepsilon$ small enough, the last
term is also smaller than $\gamma\left(  x\right)  .$ So%
\[
y\in\mathbf{B}\left(  v_{\varepsilon},\tau^{-1}\left(  \tau\left(
T_{M}\left(  y,F(x)\right)  +\varepsilon\right)  \right)  \right)  \cap
V\cap\left[  v_{\varepsilon}-\operatorname*{cone}M\right]
\]
and we can apply $\left(  iii\right)  $ for the pair $\left(  x,v_{\varepsilon
}\right)  \in\left(  U\times\mathbb{Y}\right)  \cap\operatorname*{Gr}F.$
Consequently,
\[
y\in F\left(  \mathbf{B}\left(  x,T_{M}\left(  y,F(x)\right)  +\varepsilon
\right)  \cap\left[  x+\operatorname*{cone}L\right]  \right)  .
\]
If follows that there is $x_{\varepsilon}\in B\left(  x,T_{M}\left(
y,F(x)\right)  +\varepsilon\right)  \cap\left[  x+\operatorname*{cone}%
L\right]  \cap F^{-1}\left(  y\right)  .$ We get that
\[
T_{L}\left(  x,F^{-1}\left(  y\right)  \right)  \leq T_{L}\left(
x,x_{\varepsilon}\right)  <T_{M}\left(  y,F(x)\right)  +\varepsilon.
\]
Letting $\varepsilon\rightarrow0$ we obtain $T_{L}\left(  x,F^{-1}\left(
y\right)  \right)  \leq T_{M}\left(  y,F(x)\right)  ,$ and this is the
conclusion..\hfill$\square$ 

\bigskip

We now introduce the main regularity notions used throughout the paper. Unlike
classical regularity properties, these following concepts are formulated along
successive approximation trajectories generated by admissible directions.
%We consider directional orbital regularity notions in the case $\mathbb{X}%
%=\mathbb{Y}$.
These significantly weaker notions are consistent with the framework of
\cite[pp.~356--357]{I3}. 

\begin{df}
Let $F:X\rightrightarrows X$, $\emptyset\neq L\subset S_{\mathbb{X}}$,
$U\subset X$, and $\tau>0$. The mapping $F$ is: 

\begin{enumerate}
\item[(i) ] directionally orbitally regular on $U$ w.r.t.\ $L$ with constant
$\tau$ if for every $x\in U$
\[
\tau T_{L}(x,F(x))<\mathbf{d}(x,\mathbb{X}\setminus U)\ \Longrightarrow
\ T_{L}(x,F^{-1}(x))\leq\tau T_{L}(x,F(x)).
\]

\item[(ii)] directionally orbitally Aubin continuous on $U$ w.r.t.\ $L$ with
constant $\tau$ if for every $x\in U$ and $u\in X$
\[
x\in F(u),\ \tau T_{L}(x,u)<\mathbf{d}(x,\mathbb{X}\setminus
U)\ \Longrightarrow\ T_{L}(x,F(x))\leq\tau T_{L}(x,u).
\]

\item[ (iii)] directionally orbitally open on $U$ w.r.t.\ $L$ with constant
$\tau$ if for every $x\in U$
\[
T_{L}(x,F(x))<\tau t,\ t<\mathbf{d}(x,\mathbb{X}\setminus U)\ \Longrightarrow
\ x\in F\bigl(\mathbf{B}(x,t)\cap\lbrack x+\operatorname{cone}L]\bigr).
\]

\end{enumerate}
\end{df}

\begin{pr}
\label{equiv_orb}Let $F:X\rightrightarrows X$ be a set-valued mapping. Let
$L\subset S_{\mathbb{X}}$, $U\subset X$ be nonempty sets. Take $\tau>0.\ $Then
the following statements are equivalent: 

\begin{enumerate}
\item[(i)] $F$ is directionally orbitally regular on $U$ w.r.t.\ $L$ with
constant $\tau;$ 

\item[(ii)] $F^{-1}$ is directionally orbitally Aubin continuous on $U$
w.r.t.\ $L$ with constant $\tau;$ 

\item[(iii)] $F$ is directionally orbitally open on $U$ w.r.t.\ $L$ with
constant $\tau^{-1}.$ 
\end{enumerate}
\end{pr}

\noindent Proof. $(i)$ $\Longrightarrow$ $(ii)$. Immediate since $u\in F(x)$
implies $T_{L}(x,F(x))\leq T_{L}(x,u)$. 

\smallskip$(ii)$ $\ \Longrightarrow\ $ $(iii)$. Let $x\in U$,
$t<d(x,X\setminus U)$, and $T_{L}(x,F(x))<\tau^{-1}t$. Then for any
$\varepsilon>0$ small enough there is $u_{\varepsilon}\in F\left(  x\right)
\cap\left[  x+\operatorname*{cone}L\right]  $ such that one can write
\[
\tau\left\Vert u_{\varepsilon}-x\right\Vert =\tau T_{L}\left(
x,u_{\varepsilon}\right)  <\tau\left(  T_{L}\left(  x,F\left(  x\right)
\right)  +\varepsilon\right)  <t<\mathbf{d}\left(  x,X\setminus U\right)  .
\]
By assumption, $T_{L}\left(  x,F^{-1}(x)\right)  \leq\tau T_{L}%
(x,u_{\varepsilon})<t$ whence there is $x_{\varepsilon}\in F^{-1}%
(x)\cap\left[  x+\operatorname*{cone}L\right]  $ such that $\left\Vert
x_{\varepsilon}-x\right\Vert <t.$ We deduce that $x\in F\left(  x_{\varepsilon
}\right)  \subset F\left(  \mathbf{B}\left(  x,t\right)  \cap\left[
x+\operatorname*{cone}L\right]  \right)  ,$ i.e., $\left(  iii\right)  .$ 

\smallskip$(iii)$ $\ \Longrightarrow\ $ $(i)$. Consider $x\in U$ satisfying
$\tau T_{L}(x,F(x))<d\left(  x,X\setminus U\right)  .$ Then for all
$\varepsilon>0$ small enough $T_{L}(x,F(x))<T_{L}(x,F(x))+\varepsilon
<\tau^{-1}d\left(  x,X\setminus U\right)  .$ Taking $t=\tau\left(
T_{L}(x,F(x))+\varepsilon\right)  $ we get, by our assumption, that $x\in
F\left(  \mathbf{B}\left(  x,t\right)  \cap\left[  x+\operatorname*{cone}%
L\right]  \right)  .$ Consequently, there is $u_{\varepsilon}\in\left(
\mathbf{B}\left(  x,t\right)  \cap\left[  x+\operatorname*{cone}L\right]
\right)  \cap F^{-1}\left(  x\right)  .$ We deduce that
\[
T_{L}\left(  x,F^{-1}\left(  x\right)  \right)  \leq T_{L}\left(
x,u_{\varepsilon}\right)  <t=\tau\left(  T_{L}(x,F(x))+\varepsilon\right)  .
\]
Letting $\varepsilon\rightarrow0$ we get that $(i)$ is true. The proof is
complete.\hfill$\square$
%.\hfill$\square$

\bigskip

\section{Fixed Points via Directional Regularities}

The following theorem provides existence conditions for approximate and exact
fixed points of a set-valued mapping under relaxed directional regularity
assumptions. In normed spaces, the result recovers and extends several
classical fixed point theorems from the literature, including \cite{A1},
\cite[Theorem 4]{ABB}, \cite[Theorem 5E.2]{DR-Book}, \cite[Theorems 3.1 and
3.2]{T}, as well as the celebrated Nadler fixed point theorem \cite{N}. 

\begin{thm}
\label{Main-Theorem} Let $(X,\Vert\cdot\Vert)$ be a normed space and let
$F:X\rightrightarrows X$ be a set-valued mapping. Let $U\subset X$ be a
nonempty open set and let $L\subset S_{\mathbb{X}}$ be nonempty. Assume that: 

\begin{enumerate}
\item[(i)] $F$ is directionally orbitally regular on $U$ w.r.t. $L$ and $-L$
with constant $\tau\in(0,1)$; 

\item[(ii)] there exists $\overline{x}\in U$ such that one of the following
conditions holds:
\begin{align*}
(\alpha)\qquad &  T_{L}\bigl(\overline{x},F^{-1}(\overline{x})\bigr)<(1-\tau
)\mathbf{d}\bigl(\overline{x},\mathbb{X}\setminus U\bigr),\\
(\beta)\qquad &  T_{-L}\bigl(\overline{x},F^{-1}(\overline{x})\bigr)<(1-\tau
)\mathbf{d}\bigl(\overline{x},\mathbb{X}\setminus U\bigr).
\end{align*}

\end{enumerate}

Then the following assertions hold. 

\begin{enumerate}
\item[(a)] There exists a sequence $\left(  x_{k}\right)  _{k\in\mathbb{N}%
}\subset X$ of successive approximations of $F^{-1}$ starting from
$\overline{x}$ w.r.t. $L$ such that, for every $\varepsilon>0$, there exists
$k_{\varepsilon}\in N$ satisfying
\[
x_{k}\in\operatorname*{Fix}(F)(\varepsilon)\cap U,\qquad\forall k\geq
k_{\varepsilon},
\]
and, moreover, either
\[
T_{L}\bigl(x_{k},F(x_{k})\bigr)<\varepsilon
\]
or
\[
T_{-L}\bigl(x_{k},F(x_{k})\bigr)<\varepsilon.
\]
Furthermore, if condition $(\alpha)$ holds, then
\[
\mathbf{d}\bigl(\overline{x},\operatorname*{Fix}(F)(\varepsilon)\cap
U\bigr)\leq\frac{T_{L}\bigl(\overline{x},F^{-1}(\overline{x})\bigr)}{1-\tau},
\]
whereas if condition $(\beta)$ holds, then
\[
\mathbf{d}\bigl(\overline{x},\operatorname*{Fix}(F)(\varepsilon)\cap
U\bigr)\leq\frac{T_{-L}\bigl(\overline{x},F^{-1}(\overline{x})\bigr)}{1-\tau}.
\]

\item[(b)] If $X$ is a Banach space and $F$ has closed graph, then
\[
\operatorname*{Fix}(F)\cap U\neq\emptyset.
\]
Moreover, if condition $(\alpha)$ holds, then
\[
\mathbf{d}\bigl(\overline{x},\operatorname*{Fix}(F)\cap U\bigr)\leq\frac
{T_{L}\bigl(\overline{x},F^{-1}(\overline{x})\bigr)}{1-\tau},
\]
while if condition $(\beta)$ holds, then
\[
\mathbf{d}\bigl(\overline{x},\operatorname*{Fix}(F)\cap U\bigr)\leq
\frac{T_{-L}\bigl(\overline{x},F^{-1}(\overline{x})\bigr)}{1-\tau}.
\]

\end{enumerate}
\end{thm}

\noindent Proof.  Set $x_{0}:=\overline{x}$ and assume that condition
$(\alpha)$ holds. Choose $\ell>0$ such that
\begin{equation}
T_{L}\bigl(\overline{x},F^{-1}(\overline{x})\bigr)<\ell(1-\tau)<(1-\tau
)\mathbf{d}\bigl(\overline{x},\mathbb{X}\setminus U\bigr).\label{OK123}%
\end{equation}

By the definition of $T_{L}$, there exist $t_{0}\geq0$ and $u_{0}\in L$ such
that
\[
x_{1}:=\overline{x}+t_{0}u_{0}\in F^{-1}(\overline{x})
\]
and
\[
t_{0}<\ell(1-\tau).
\]
Equivalently,
\[
x_{0}=\overline{x}\in F(x_{1}),
\]
and
\[
\Vert x_{1}-x_{0}\Vert=t_{0}<(1-\tau)\mathbf{d}\bigl(\overline{x}%
,\mathbb{X}\setminus U\bigr).
\]
Since $U$ is open, it follows that $x_{1}\in U$. Moreover,
\[
x_{0}\in F(x_{1})\cap\bigl(x_{1}-\operatorname{cone}L\bigr),
\]
and therefore
\[
T_{-L}\bigl(x_{1},F(x_{1})\bigr)\leq\Vert x_{1}-x_{0}\Vert.
\]

Next, using the fact that the distance function $d(\cdot,X\setminus U)$ is
$1$-Lipschitz, we obtain
\[
\Vert x_{1}-x_{0}\Vert<(1-\tau)\mathbf{d}\bigl(x_{1},\mathbb{X}\setminus
U\bigr)+(1-\tau)\Vert x_{1}-\overline{x}\Vert.
\]
Consequently,
\[
\tau\Vert x_{1}-x_{0}\Vert<(1-\tau)\mathbf{d}\bigl(x_{1},\mathbb{X}\setminus
U\bigr).
\]
Combining this estimate with the previous inequality yields
\[
\tau T_{-L}\bigl(x_{1},F(x_{1})\bigr)<(1-\tau)\mathbf{d}\bigl(x_{1}%
,\mathbb{X}\setminus U\bigr)<\mathbf{d}\bigl(x_{1},\mathbb{X}\setminus
U\bigr).
\]

Applying assumption (i) for $-L$, we deduce that
\[
T_{-L}\bigl(x_{1},F^{-1}(x_{1})\bigr)\leq\tau T_{-L}\bigl(x_{1},F(x_{1}%
)\bigr).
\]
Hence,
\[
T_{-L}\bigl(x_{1},F^{-1}(x_{1})\bigr)<\min\left\{  \ell\tau(1-\tau
),(1-\tau)\mathbf{d}\bigl(x_{1},\mathbb{X}\setminus U\bigr)\right\}  .
\]

Therefore, there exist $t_{1}\geq0$ and $u_{1}\in L$ such that
\[
x_{2}:=x_{1}-t_{1}u_{1}\in F^{-1}(x_{1}),
\]
with
\[
\Vert x_{2}-x_{1}\Vert=t_{1}<\min\left\{  \ell\tau(1-\tau),(1-\tau
)\mathbf{d}\bigl(x_{1},\mathbb{X}\setminus U\bigr)\right\}  .
\]
Again, this implies that $x_{2}\in U$, and
\[
T_{L}\bigl(x_{2},F(x_{2})\bigr)\leq\Vert x_{2}-x_{1}\Vert.
\]
Proceeding exactly as above, we obtain
\[
\tau\Vert x_{2}-x_{1}\Vert<(1-\tau)\mathbf{d}\bigl(x_{2},\mathbb{X}\setminus
U\bigr),
\]
and therefore, by assumption (i) applied to $L$,
\[
T_{L}\bigl(x_{2},F^{-1}(x_{2})\bigr)\leq\tau T_{L}\bigl(x_{2},F(x_{2}%
)\bigr)<\min\left\{  \ell\tau^{2}(1-\tau),(1-\tau)\mathbf{d}\bigl(x_{2}%
,\mathbb{X}\setminus U\bigr)\right\}  .
\]

Proceeding inductively, we construct sequences $\left(  t_{n}\right)  \subset
R_{+}$, $\left(  u_{n}\right)  \subset L\cup(-L)$, and $\left(  x_{n}\right)
\subset X$ such that, for every $n\in N$,
\begin{align}
&  x_{n+1}=x_{n}+t_{n}u_{n},\qquad x_{n}\in F(x_{n+1})\cap U,\label{iter1}\\
&  0\leq t_{n}=\Vert x_{n+1}-x_{n}\Vert<\min\left\{  \ell\tau^{n}%
(1-\tau),(1-\tau)\mathbf{d}(x_{n},\mathbb{X}\setminus U)\right\}
,\label{iter2}\\
&  u_{n}\in-L\quad\text{for odd }n,\qquad u_{n}\in L\quad\text{for even
}n.\label{iter3}%
\end{align}

Relations \eqref{iter1}--\eqref{iter3} imply that
\[
\mathbf{d}\bigl(x_{n+1},F(x_{n+1})\bigr)\leq\Vert x_{n+1}-x_{n}\Vert
=t_{n}<\ell\tau^{n}(1-\tau).
\]
Moreover, since
\[
x_{n+1}+t_{n}(-u_{n})=x_{n}\in F(x_{n+1}),
\]
we have:
\[
T_{L}\bigl(x_{n+1},F(x_{n+1})\bigr)\leq t_{n}<\ell\tau^{n}(1-\tau
)\qquad\text{for odd }n,
\]
and
\[
T_{-L}\bigl(x_{n+1},F(x_{n+1})\bigr)\leq t_{n}<\ell\tau^{n}(1-\tau
)\qquad\text{for even }n.
\]

Consequently, for every $\varepsilon>0$, there exists $n_{\varepsilon}\in N$
such that, for all $n\geq n_{\varepsilon}$,
\[
x_{n+1}\in\operatorname*{Fix}(F)(\varepsilon)\cap U,
\]
and either
\[
T_{L}\bigl(x_{n+1},F(x_{n+1})\bigr)<\varepsilon
\]
or
\[
T_{-L}\bigl(x_{n+1},F(x_{n+1})\bigr)<\varepsilon.
\]

Next, for $m>n$,
\begin{align*}
\Vert x_{m}-x_{n}\Vert &  \leq\sum_{k=n}^{m-1}\Vert x_{k+1}-x_{k}\Vert\\
&  <\sum_{k=n}^{m-1}\ell\tau^{k}(1-\tau)\\
&  =\ell\tau^{n}(1-\tau)\sum_{k=0}^{m-n-1}\tau^{k}\\
&  =\ell\tau^{n}(1-\tau^{m-n}).
\end{align*}
Thus, $\left(  x_{n}\right)  $ is a Cauchy sequence in $X$. In particular,
taking $n=0$, we obtain
\begin{equation}
\Vert\overline{x}-x_{m}\Vert<\ell(1-\tau^{m}).\label{OK9}%
\end{equation}
Therefore,
\[
\mathbf{d}\bigl(\overline{x},\operatorname*{Fix}(F)(\varepsilon)\cap
U\bigr)\leq\Vert\overline{x}-x_{m}\Vert<\ell.
\]
Since $\ell$ can be chosen arbitrarily close to
\[
\frac{T_{L}\bigl(\overline{x},F^{-1}(\overline{x})\bigr)}{1-\tau},
\]
it follows that
\[
\mathbf{d}\bigl(\overline{x},\operatorname*{Fix}(F)(\varepsilon)\cap
U\bigr)\leq\frac{T_{L}\bigl(\overline{x},F^{-1}(\overline{x})\bigr)}{1-\tau}.
\]

Assume now that $X$ is complete and that $\operatorname{Gr}F$ is closed. Since
$\left(  x_{n}\right)  $ is Cauchy, there exists $\hat{x}\in X$ such that
$x_{n}\rightarrow\hat{x}$. Because
\[
(x_{n+1},x_{n})\rightarrow(\hat{x},\hat{x})
\]
and
\[
x_{n}\in F(x_{n+1}),
\]
the closedness of $\operatorname{Gr}F$ yields
\[
\hat{x}\in F(\hat{x}),
\]
that is,
\[
\hat{x}\in\operatorname*{Fix}(F).
\]

Passing to the limit in \eqref{OK9}, we obtain
\[
\Vert\overline{x}-\hat{x}\Vert\leq\ell<\mathbf{d}\bigl(\overline{x}%
,\mathbb{X}\setminus U\bigr).
\]
Since $U$ is open, this implies that $\hat{x}\in U$. Therefore,
\[
\mathbf{d}\bigl(\overline{x},\operatorname*{Fix}(F)\cap U\bigr)\leq
\Vert\overline{x}-\hat{x}\Vert\leq\ell.
\]
Letting
\[
\ell\downarrow\frac{T_{L}\bigl(\overline{x},F^{-1}(\overline{x})\bigr)}%
{1-\tau},
\]
we conclude that
\[
\mathbf{d}\bigl(\overline{x},\operatorname*{Fix}(F)\cap U\bigr)\leq\frac
{T_{L}\bigl(\overline{x},F^{-1}(\overline{x})\bigr)}{1-\tau}.
\]

This completes the proof under assumption $(\alpha)$. The proof under
assumption $(\beta)$ is entirely analogous, with the roles of even and odd
indices interchanged in the construction of the sequence $\left(
x_{n}\right)  $. \hfill$\square$ 

\bigskip

Building upon the previous results and their proofs, together with
Proposition~\ref{equiv_orb}, we obtain the following theorem. 

\begin{thm}
\label{Direct Lip} Let $\left(  \mathbb{X},\Vert\cdot\Vert\right)  $ be a
normed space, and let $G:X\rightrightarrows X$ be a set-valued mapping. Let
$U\subset X$ be a nonempty open set and let $L\subset S_{\mathbb{X}}$ be
nonempty. Assume that: 

\begin{enumerate}
\item[(i)] $G$ is directionally orbitally Aubin continuous on $U$ w.r.t. $L$
and $-L$ with constant $\tau\in(0,1)$; 

\item[(ii)] there exists $\overline{x}\in U$ such that one of the following
conditions holds:
\begin{align*}
(\alpha)\qquad &  T_{L}\bigl(\overline{x},G(\overline{x})\bigr)<(1-\tau
)\mathbf{d}\bigl(\overline{x},\mathbb{X}\setminus U\bigr),\\
(\beta)\qquad &  T_{-L}\bigl(\overline{x},G(\overline{x})\bigr)<(1-\tau
)\mathbf{d}\bigl(\overline{x},\mathbb{X}\setminus U\bigr).
\end{align*}

\end{enumerate}

Then the following assertions hold. 

\begin{enumerate}
\item[(a)] There exists a sequence $\left(  x_{k}\right)  _{k\in\mathbb{N}%
}\subset X$ of successive approximations of $G$ starting from $\overline{x}$
w.r.t. $L$ such that, for every $\varepsilon>0$, there exists $k_{\varepsilon
}\in N$ satisfying
\[
x_{k}\in\operatorname*{Fix}(G)(\varepsilon)\cap U,\qquad\forall k\geq
k_{\varepsilon},
\]
and, moreover, either
\[
T_{L}\bigl(x_{k},G(x_{k})\bigr)\leq\varepsilon
\]
or
\[
T_{-L}\bigl(x_{k},G(x_{k})\bigr)\leq\varepsilon.
\]
Furthermore, if condition $(\alpha)$ holds, then
\[
\mathbf{d}\bigl(\overline{x},\operatorname*{Fix}(G)(\varepsilon)\cap
U\bigr)\leq\frac{T_{L}\bigl(\overline{x},G(\overline{x})\bigr)}{1-\tau},
\]
whereas if condition $(\beta)$ holds, then
\[
\mathbf{d}\bigl(\overline{x},\operatorname*{Fix}(G)(\varepsilon)\cap
U\bigr)\leq\frac{T_{-L}\bigl(\overline{x},G(\overline{x})\bigr)}{1-\tau}.
\]

\item[(b)] If $X$ is a Banach space and $G$ has closed graph, then
\[
\operatorname*{Fix}(G)\cap U\neq\emptyset.
\]
Moreover, if condition $(\alpha)$ holds, then
\[
\mathbf{d}\bigl(\overline{x},\operatorname*{Fix}(G)\cap U\bigr)\leq\frac
{T_{L}\bigl(\overline{x},G(\overline{x})\bigr)}{1-\tau},
\]
while if condition $(\beta)$ holds, then
\[
\mathbf{d}\bigl(\overline{x},\operatorname*{Fix}(G)\cap U\bigr)\leq
\frac{T_{-L}\bigl(\overline{x},G(\overline{x})\bigr)}{1-\tau}.
\]

\end{enumerate}
\end{thm}

\noindent Proof.  By Proposition~\ref{equiv_orb}, the inverse mapping $G^{-1}$
is directionally orbitally regular on $U$ w.r.t. $L$ and $-L$ with the same
constant $\tau$. 

Applying Theorem~\ref{Main-Theorem} to the mapping $G^{-1}$, we obtain a
sequence of successive approximations $\left(  x_{n}\right)  _{n\in\mathbb{N}%
}$ satisfying
\[
x_{n+1}=x_{n}+t_{n}u_{n},\qquad x_{n}\in G^{-1}(x_{n+1}),
\]
where $u_{n}\in\pm L$ and $t_{n}\geq0$. 

Since
\[
x_{n}\in G^{-1}(x_{n+1})\quad\Longleftrightarrow\quad x_{n+1}\in G(x_{n}),
\]
it follows that $\left(  x_{n}\right)  $ is a sequence of successive
approximations of $G$. Moreover,
\[
x_{n}=x_{n+1}+t_{n}(-u_{n}),
\]
and therefore
\[
T_{\mp L}\bigl(x_{n},G(x_{n})\bigr)\leq t_{n}.
\]

The remaining conclusions follow exactly as in the proof of
Theorem~\ref{Main-Theorem}. \hfill$\square$ 

\begin{rmk}
\label{remark12} Observe that, in Theorem~\ref{Main-Theorem}, if
$\operatorname{cone}L$ is convex, then the elements of $\operatorname*{Fix}%
(F)(\varepsilon)\cap U$ obtained in part~(a) necessarily belong to
\[
\overline{x}+\operatorname{cone}L-\operatorname{cone}L.
\]

\end{rmk}

\begin{rmk}
\label{rem_eta} A closer inspection of the proof of Theorem~\ref{Main-Theorem}
shows that the sequence of successive approximations $\left(  x_{k}\right)
_{k\in\mathbb{N}}$ may be chosen so as to satisfy
\[
\Vert x_{n}-\overline{x}\Vert\leq\frac{T_{L}\bigl(\overline{x},F^{-1}%
(\overline{x})\bigr)+\eta}{1-\tau},\qquad\forall n\in\mathbb{N},
\]
for every fixed constant $\eta>0$. 

Indeed, it suffices to choose $\rho>0$ such that
\[
T_{L}\bigl(\overline{x},F^{-1}(\overline{x})\bigr)+\rho<(1-\tau)\mathbf{d}%
\bigl(\overline{x},\mathbb{X}\setminus U\bigr),
\]
and then select $\ell>0$ satisfying
\[
T_{L}\bigl(\overline{x},F^{-1}(\overline{x})\bigr)+\rho<\ell(1-\tau
)<(1-\tau)\mathbf{d}\bigl(\overline{x},\mathbb{X}\setminus U\bigr).
\]

An analogous observation applies to the sequence of successive approximations
constructed in Theorem~\ref{Direct Lip}. Namely, for every $\eta>0$, the
sequence may be chosen so that
\[
\Vert x_{n+1}-\overline{x}\Vert\leq\frac{T_{L}\bigl(\overline{x}%
,G(\overline{x})\bigr)+\eta}{1-\tau},\qquad\forall n\in\mathbb{N}.
\]

In particular, these remarks show that our results extend
\cite[Proposition~5.1]{ABB}. 
\end{rmk}

By taking
\[
U:=\mathbf{B}(\overline{x},r),\qquad r>0,
\]
in Theorems~\ref{Main-Theorem} and~\ref{Direct Lip}, one immediately obtains
local versions of these results. 

In the global setting, namely when $U=X$, additional care is required, since
the minimal time function may attain the value $+\infty$. In particular,
condition~(ii) in the previous theorems may hold only for specific points. To
make these modifications explicit, we state below a global version of
Theorem~\ref{Direct Lip}. 

\begin{thm}
\label{Direct Lip-Global} Let $\left(  \mathbb{X},\Vert\cdot\Vert\right)  $ be
a normed space, let $G:X\rightrightarrows X$ be a set-valued mapping, and let
$L\subset S_{\mathbb{X}}$ be nonempty. Assume that: 

\begin{enumerate}
\item[(i)] $G$ is globally directionally orbital Aubin continuous w.r.t. $L$
and $-L$ with constant $\tau\in(0,1)$; 

\item[(ii)] there exists $\overline{x}\in X$ such that
\[
\overline{x}\in G(\overline{x})+\operatorname{cone}(L\cup(-L)).
\]

\end{enumerate}

Then all conclusions of Theorem~\ref{Direct Lip} remain valid. 
\end{thm}

\noindent Proof.  When $U=X$, condition~(ii) in Theorem~\ref{Direct Lip}
reduces to the requirement that either
\[
\overline{x}\in G(\overline{x})+\operatorname{cone}L
\]
or
\[
\overline{x}\in G(\overline{x})+\operatorname{cone}(-L).
\]
Equivalently,
\[
\overline{x}\in G(\overline{x})+\operatorname{cone}(L\cup(-L)),
\]
which is precisely assumption~(ii). 

The proof then follows the same lines as that of Theorem~\ref{Direct Lip}. In
particular, one chooses
\[
\ell>T_{L}\bigl(\overline{x},G(\overline{x})\bigr)
\]
arbitrarily and repeats the iterative construction. \hfill$\square$ 

\begin{rmk}
Condition~(ii) in Theorem~\ref{Direct Lip-Global} is automatically satisfied
whenever $L=S_{\mathbb{X}}$ and the mapping $G$ has nonempty values. 
\end{rmk}

\begin{rmk}
[Comparison with existing results]When $L=S_{\mathbb{X}}$ and the minimal time
function reduces to the metric distance, our results recover several classical
fixed point theorems, including Nadler's theorem and variants of the
Lyusternik--Graves principle. In contrast with contraction-type approaches,
our assumptions are directional and orbital in nature and do not require
global Lipschitz behavior. 
\end{rmk}

\section{Applications}

In this section we discuss several consequences of the preceding developments,
focusing on coincidence theorems and stability properties for coupled
set-valued systems. 

\subsection{Applications to coincidence results}

The aim of this subsection is to derive coincidence, coupled fixed point, and
stability results from the directional regularity framework developed above in
the setting of normed spaces. 

We begin with a simple observation describing the behavior of minimal time
functions under a rescaling of the ambient norm. Let $\left(  \mathbb{X}%
,\Vert\cdot\Vert\right)  $ be a normed space, $L\subset S_{\mathbb{X}}$, and
let $\alpha>0$. Endow $X$ with the equivalent norm $\alpha\Vert\cdot\Vert$.
Then $\alpha^{-1}L\subset S_{(\mathbb{X},\alpha\Vert\cdot\Vert)}$. For every
nonempty set $A\subset X$, the associated directional minimal time functions
satisfy
\[
T_{\alpha^{-1}L}^{\,\alpha\Vert\cdot\Vert}(\cdot,A)=\alpha\,T_{L}%
^{\,\Vert\cdot\Vert}(\cdot,A).
\]

Let $X$ and $Y$ be normed spaces, and let $\tau_{1},\tau_{2}>0$. Throughout
this subsection we consider the product space $X\times Y$ equipped with the
norms $\Vert\cdot\Vert_{1},\Vert\cdot\Vert_{2},\Vert\cdot\Vert_{3}$ defined
by
\[
\Vert(x,y)\Vert_{1}:=\max\{\Vert x\Vert,\Vert y\Vert\},
\]%
\[
\Vert(x,y)\Vert_{2}:=\max\left\{  \sqrt{\tau_{1}^{-1}}\Vert x\Vert,\sqrt
{\tau_{2}^{-1}}\Vert y\Vert\right\}  ,
\]
and
\[
\Vert(x,y)\Vert_{3}:=\max\left\{  \sqrt{\tau_{1}}\Vert x\Vert,\sqrt{\tau_{2}%
}\Vert y\Vert\right\}  .
\]

We shall work under the following standing assumptions. 

\medskip

\noindent(S)\  $X$ and $Y$ are normed spaces, $F_{1}:X\rightrightarrows Y$ and
$F_{2}:Y\rightrightarrows X$ are set-valued mappings, $L\subset S_{\mathbb{X}%
}$ and $M\subset S_{\mathbb{Y}}$ are nonempty sets of directions, and
$F:X\times Y\rightrightarrows X\times Y$ is defined by
\[
\mathcal{F}(x,y):=(F_{2}(y),F_{1}(x)).
\]

\begin{pr}
\label{pr1} Assume that the setting (S) holds. Let $U_{1},U_{2}\subset X$ and
$V_{1},V_{2}\subset Y$ be open sets. Let $\gamma$ and $\delta$ be gauge
functions on $U_{1}$ and $V_{2}$, respectively. Suppose that there exist
$\tau_{1},\tau_{2}>0$ such that: 

\begin{itemize}
\item[(i)] $F_{1}$ is directionally $\gamma$-regular on $U_{1}\times V_{1}$
w.r.t. $L$ and $M$ with constant $\tau_{1}$; 

\item[(ii)] $F_{2}$ is directionally $\delta$-regular on $V_{2}\times U_{2}$
w.r.t. $M$ and $L$ with constant $\tau_{2}$. 
\end{itemize}

Then the following assertions hold. 

\begin{itemize}
\item[(a)] On the space $\bigl(X\times Y,\Vert\cdot\Vert_{1}\bigr)$, the
mapping $F$ is directionally $\rho_{1}$-regular on $(U_{1}\times V_{2}%
)\times(U_{2}\times V_{1})$ w.r.t. $L\times M$ and $L\times M$ with constant
\[
\tau:=\max\{\tau_{1},\tau_{2}\},
\]
where
\[
\rho_{1}(x,y):=\min\{\gamma(x),\delta(y)\}.
\]

\item[(b)] On the space $\bigl(X\times Y,\Vert\cdot\Vert_{2}\bigr)$, the
mapping $F$ is directionally $\rho_{2}$-regular on $(U_{1}\times V_{2}%
)\times(U_{2}\times V_{1})$ w.r.t. $\sqrt{\tau_{1}}L\times\sqrt{\tau_{2}}M$
and $\sqrt{\tau_{1}}L\times\sqrt{\tau_{2}}M$ with constant $\sqrt{\tau_{1}%
\tau_{2}}$, where
\[
\rho_{2}(x,y):=\min\left\{  \sqrt{\tau_{1}^{-1}}\gamma(x),\sqrt{\tau_{2}^{-1}%
}\delta(y)\right\}  .
\]

\end{itemize}
\end{pr}

\noindent Proof.  (a) Let $(x,y)\in U_{1}\times V_{2}$ and $(u,v)\in
U_{2}\times V_{1}$ satisfy
\[
\tau T_{L\times M}\bigl((u,v),\mathcal{F}(x,y)\bigr)<\rho_{1}(x,y).
\]
Equivalently,
\[
\max\{\tau_{1},\tau_{2}\}\max\left\{  T_{L}\bigl(u,F_{2}(y)\bigr),T_{M}%
\bigl(v,F_{1}(x)\bigr)\right\}  <\min\{\gamma(x),\delta(y)\}.
\]
Hence,
\[
\tau_{2}T_{L}\bigl(u,F_{2}(y)\bigr)<\delta(y)
\]
and
\[
\tau_{1}T_{M}\bigl(v,F_{1}(x)\bigr)<\gamma(x).
\]
By assumptions (i) and (ii),
\[
T_{L}\bigl(x,F_{1}^{-1}(v)\bigr)\leq\tau_{1}T_{M}\bigl(v,F_{1}(x)\bigr)
\]
and
\[
T_{M}\bigl(y,F_{2}^{-1}(u)\bigr)\leq\tau_{2}T_{L}\bigl(u,F_{2}(y)\bigr).
\]
Therefore,
\begin{align*}
T_{L\times M}\bigl((x,y),\mathcal{F}^{-1}(u,v)\bigr) &  =\max\left\{
T_{L}\bigl(x,F_{1}^{-1}(v)\bigr),T_{M}\bigl(y,F_{2}^{-1}(u)\bigr)\right\}  \\
&  \leq\max\{\tau_{1},\tau_{2}\}\max\left\{  T_{L}\bigl(u,F_{2}(y)\bigr),T_{M}%
\bigl(v,F_{1}(x)\bigr)\right\}  \\
&  =\tau T_{L\times M}\bigl((u,v),\mathcal{F}(x,y)\bigr).
\end{align*}
This proves assertion (a). 

\medskip

(b) Consider now the space $\bigl(X\times Y,\Vert\cdot\Vert_{2}\bigr)$. Let
$(x,y)\in U_{1}\times V_{2}$ and $(u,v)\in U_{2}\times V_{1}$ satisfy
\[
\sqrt{\tau_{1}\tau_{2}}\,T_{\sqrt{\tau_{1}}L\times\sqrt{\tau_{2}}M}%
^{\,\Vert\cdot\Vert_{2}}\bigl((u,v),\mathcal{F}(x,y)\bigr)<\rho_{2}(x,y).
\]
Using the scaling property of the minimal time function, this inequality can
be rewritten as
\[
\sqrt{\tau_{1}\tau_{2}}\max\left\{  \sqrt{\tau_{1}^{-1}}T_{L}\bigl(u,F_{2}%
(y)\bigr),\sqrt{\tau_{2}^{-1}}T_{M}\bigl(v,F_{1}(x)\bigr)\right\}
<\min\left\{  \sqrt{\tau_{1}^{-1}}\gamma(x),\sqrt{\tau_{2}^{-1}}%
\delta(y)\right\}  .
\]
Consequently,
\[
\tau_{2}T_{L}\bigl(u,F_{2}(y)\bigr)<\delta(y)
\]
and
\[
\tau_{1}T_{M}\bigl(v,F_{1}(x)\bigr)<\gamma(x).
\]
Applying assumptions (i) and (ii), we obtain
\[
T_{L}\bigl(x,F_{1}^{-1}(v)\bigr)\leq\tau_{1}T_{M}\bigl(v,F_{1}(x)\bigr)
\]
and
\[
T_{M}\bigl(y,F_{2}^{-1}(u)\bigr)\leq\tau_{2}T_{L}\bigl(u,F_{2}(y)\bigr).
\]
Hence,
\begin{align*}
&  T_{\sqrt{\tau_{1}}L\times\sqrt{\tau_{2}}M}^{\,\Vert\cdot\Vert_{2}%
}\bigl((x,y),\mathcal{F}^{-1}(u,v)\bigr)\\
&  =\max\left\{  \sqrt{\tau_{1}^{-1}}T_{L}\bigl(x,F_{1}^{-1}(v)\bigr),\sqrt
{\tau_{2}^{-1}}T_{M}\bigl(y,F_{2}^{-1}(u)\bigr)\right\}  \\
&  \leq\max\left\{  \sqrt{\tau_{1}}T_{M}\bigl(v,F_{1}(x)\bigr),\sqrt{\tau_{2}%
}T_{L}\bigl(u,F_{2}(y)\bigr)\right\}  \\
&  \leq\sqrt{\tau_{1}\tau_{2}}\max\left\{  \sqrt{\tau_{1}^{-1}}T_{L}%
\bigl(u,F_{2}(y)\bigr),\sqrt{\tau_{2}^{-1}}T_{M}\bigl(v,F_{1}%
(x)\bigr)\right\}  \\
&  =\sqrt{\tau_{1}\tau_{2}}\,T_{\sqrt{\tau_{1}}L\times\sqrt{\tau_{2}}%
M}^{\,\Vert\cdot\Vert_{2}}\bigl((u,v),\mathcal{F}(x,y)\bigr).
\end{align*}
The proof is complete. \hfill$\square$ 

\begin{pr}
\label{pr2}In the setting (S), let $U,$ $V$ be open subsets in $X$ and $Y.$
Assume that there exist $\tau_{1},\tau_{2}>0$ such that 

\begin{enumerate}
\item[(i)] $F_{1}$ is Milyutin directionally regular on $U\times V$w.r.t.\ $L$
and $M$ with constant $\tau_{1}$; 

\item[(ii)] $F_{2}$ is Milyutin directionally regular on $V\times U$
w.r.t.\ $M$ and $L$ with constant $\tau_{2}$. 
\end{enumerate}

Then the following assertions hold. 

\begin{enumerate}
\item[(a)] On the space $(X\times Y,\left\Vert \cdot\right\Vert _{1}),$ $F$ is
directionally orbitally regular on $U\times V$ w.r.t.\ $L\times M$ with
constant $\tau:=\max\left\{  \tau_{1},\tau_{2}\right\}  $. 

\item[(b)] On the space $(X\times Y,\left\Vert \cdot\right\Vert _{2}),$ $F$ is
directionally orbitally regular on $U\times V$ w.r.t.\ $\sqrt{\tau_{1}}%
L\times\sqrt{\tau_{2}}M$ with constant $\sqrt{\tau_{1}\tau_{2}}.$ 
\end{enumerate}
\end{pr}

\noindent Proof. (a) Let $\left(  x,y\right)  \in U\times V$ such that
\[
\tau T_{L\times M}\left(  \left(  x,y\right)  ,\mathcal{F}\left(  x,y\right)
\right)  <\mathbf{d}\left(  \left(  x,y\right)  ,\left(  \mathbb{X}%
\times\mathbb{Y}\right)  \setminus\left(  U\times V\right)  \right)  .
\]

Since
\[
\mathbf{d}\left(  \left(  x,y\right)  ,\mathbb{X}\times\mathbb{Y}\setminus
U\times V\right)  =\min\left\{  \mathbf{d}\left(  x,\mathbb{X}\setminus
U\right)  ,\mathbf{d}\left(  y,\mathbb{Y}\setminus V\right)  \right\}  ,
\]
we get that
\[
\tau T_{L}\left(  x,F_{2}\left(  y\right)  \right)  <\mathbf{d}\left(
y,\mathbb{Y}\setminus V\right)  \text{,}%
\]
and
\[
\tau T_{M}\left(  y,F_{1}\left(  x\right)  \right)  <\mathbf{d}\left(
x,\mathbb{X}\setminus U\right)  .
\]
Applying assumptions $(i)$ and $(ii)$ we obtain%
\[
T_{L}\left(  x,F_{1}^{-1}\left(  y\right)  \right)  <\tau T_{M}\left(
y,F_{1}\left(  x\right)  \right)  ,
\]
and
\[
T_{M}\left(  y,F_{2}^{-1}\left(  x\right)  \right)  <\tau T_{L}\left(
x,F_{2}\left(  y\right)  \right)  .
\]
Therefore
\begin{align*}
T_{L\times M}\left(  \left(  x,y\right)  ,\mathcal{F}^{-1}(x,y)\right)   &
=T_{L\times M}\left(  \left(  x,y\right)  ,\left(  F_{1}^{-1}\left(  y\right)
,F_{2}\left(  x\right)  \right)  \right)  \\
&  =\max\left\{  T_{L}\left(  x,F_{1}^{-1}\left(  y\right)  \right)
,T_{M}\left(  y,F_{2}^{-1}\left(  x\right)  \right)  \right\}  \\
&  \leq\tau\max\left\{  T_{M}\left(  y,F_{1}\left(  x\right)  \right)
,T_{L}\left(  x,F_{2}\left(  y\right)  \right)  \right\}  \\
&  =\tau T_{L\times M}\left(  \left(  x,y\right)  ,\mathcal{F}\left(
x,y\right)  \right)  ,
\end{align*}
and we obtain the announced property. 

(b) We consider now the space $(X\times Y,\left\Vert \cdot\right\Vert _{2})$.
Let $(x,y)\in U\times V$ be such that
\[
\sqrt{\tau_{1}\tau_{2}}T_{\sqrt{\tau_{1}}L\times\sqrt{\tau_{2}}M}^{\left\Vert
\cdot\right\Vert _{2}}((x,y),\mathcal{F}(x,y))<\mathbf{d}\left(  \left(
x,y\right)  ,\left(  \mathbb{X}\times\mathbb{Y}\right)  \setminus\left(
U\times V\right)  \right)  .
\]
Because of
\[
\mathbf{d}\left(  \left(  x,y\right)  ,\mathbb{X}\times\mathbb{Y}\setminus
U\times V\right)  =\min\left\{  \sqrt{\tau_{1}^{-1}}\mathbf{d}\left(
x,\mathbb{X}\setminus U\right)  ,\sqrt{\tau_{2}^{-1}}\mathbf{d}\left(
y,\mathbb{Y}\setminus V\right)  \right\}  ,
\]
the above inequality becomes
\begin{align*}
&  \sqrt{\tau_{1}\tau_{2}}\max\left\{  \sqrt{\tau_{1}^{-1}}T_{L}\left(
x,F_{2}\left(  y\right)  \right)  ,\sqrt{\tau_{2}^{-1}}T_{M}\left(
y,F_{1}\left(  x\right)  \right)  \right\}  \\
&  <\min\left\{  \sqrt{\tau_{1}^{-1}}\mathbf{d}\left(  x,\mathbb{X}\setminus
U\right)  ,\sqrt{\tau_{2}^{-1}}\mathbf{d}\left(  y,\mathbb{Y}\setminus
V\right)  \right\}
\end{align*}
which implies
\[
\tau_{2}T_{L}\left(  x,F_{2}\left(  y\right)  \right)  <\mathbf{d}\left(
y,\mathbb{Y}\setminus V\right)  ,\text{ }%
\]
and
\[
\tau_{1}T_{M}\left(  y,F_{1}\left(  x\right)  \right)  <\mathbf{d}\left(
x,\mathbb{X}\setminus U\right)  .
\]
Again by $\left(  i\right)  $ and $\left(  ii\right)  ,$
\[
T_{L}\left(  x,F_{1}^{-1}\left(  y\right)  \right)  <\tau_{1}T_{M}\left(
y,F_{1}\left(  x\right)  \right)  ,\text{ }%
\]
and
\[
T_{M}\left(  y,F_{2}^{-1}\left(  x\right)  \right)  <\tau_{2}T_{L}\left(
x,F_{2}\left(  y\right)  \right)  .
\]
We get that
\begin{align*}
T_{\sqrt{\tau_{1}}L\times\sqrt{\tau_{2}}M}^{\left\Vert \cdot\right\Vert _{2}%
}\left(  \left(  x,y\right)  ,\mathcal{F}^{-1}(x,y)\right)   &  =T_{\sqrt
{\tau_{1}}L\times\sqrt{\tau_{2}}M}^{\left\Vert \cdot\right\Vert _{2}}\left(
\left(  x,y\right)  ,\left(  F_{1}^{-1}\left(  v\right)  ,F_{2}^{-1}\left(
u\right)  \right)  \right)  \\
&  =\max\left\{  \sqrt{\tau_{1}^{-1}}T_{L}\left(  x,F_{1}^{-1}\left(
y\right)  \right)  ,\sqrt{\tau_{2}^{-1}}T_{M}\left(  y,F_{2}^{-1}\left(
u\right)  \right)  \right\}  \\
&  \leq\max\left\{  \sqrt{\tau_{1}}T_{M}\left(  y,F_{1}\left(  x\right)
\right)  ,\sqrt{\tau_{2}}T_{L}\left(  x,F_{2}\left(  y\right)  \right)
\right\}  \\
&  \leq\sqrt{\tau_{1}\tau_{2}}\max\left\{  \sqrt{\tau_{1}^{-1}}T_{L}\left(
x,F_{2}\left(  y\right)  \right)  ,\sqrt{\tau_{2}^{-1}}T_{M}\left(
y,F_{1}\left(  x\right)  \right)  \right\}  \\
&  =\sqrt{\tau_{1}\tau_{2}}T_{\sqrt{\tau_{1}}L\times\sqrt{\tau_{2}}%
M}^{\left\Vert \cdot\right\Vert _{2}}((x,y),\mathcal{F}(x,y)).
\end{align*}
Therefore $F$ is directionally orbitally regular on $U\times V$ with respect
to $\sqrt{\tau_{1}}L\times\sqrt{\tau_{2}}M$ with constant $\sqrt{\tau_{1}%
\tau_{2}}.$\hfill$\square$ 

\bigskip

The next result shows that directional orbital Aubin continuity is preserved
under the product construction introduced in the setting (S). The proof
follows the same scheme as that of Proposition~\ref{pr1}, but requires a
suitable adaptation to the orbital Aubin framework. 

\begin{pr}
\label{pr3} In the setting (S), let $U\subset X$ and $V\subset Y$ be open
sets. Assume that there exist constants $\tau_{1},\tau_{2}>0$ such that: 

\begin{enumerate}
\item[(i)] $F_{1}$ is Milyutin directionally Aubin continuous on $U\times V$
w.r.t. $L$ and $M$ with constant $\tau_{1}$; 

\item[(ii)] $F_{2}$ is Milyutin directionally Aubin continuous on $V\times U$
w.r.t. $M$ and $L$ with constant $\tau_{2}$. 
\end{enumerate}

Then the following assertions hold. 

\begin{enumerate}
\item[(a)] On the product space $\left(  \mathbb{X}\times\mathbb{Y},\Vert
\cdot\Vert_{1}\right)  $, the mapping $F$ is directionally orbitally Aubin
continuous on $U\times V$ w.r.t. $L\times M$ with constant
\[
\tau:=\max\{\tau_{1},\tau_{2}\}.
\]

\item[(b)] On the product space $\left(  \mathbb{X}\times\mathbb{Y},\Vert
\cdot\Vert_{3}\right)  $, the mapping $F$ is directionally orbitally Aubin
continuous on $U\times V$ w.r.t.
\[
\sqrt{\tau_{1}^{-1}}\,L\times\sqrt{\tau_{2}^{-1}}\,M
\]
with constant $\sqrt{\tau_{1}\tau_{2}}$. 
\end{enumerate}
\end{pr}

\noindent Proof.  (a) Let $(x,y)\in U\times V$ and $(u,v)\in X\times Y$
satisfy
\[
(x,y)\in\mathcal{F}(u,v)
\]
and
\[
\max\{\tau_{1},\tau_{2}\}T_{L\times M}\bigl((x,y),(u,v)\bigr)<\mathbf{d}%
\bigl((x,y),(\mathbb{X}\times\mathbb{Y})\setminus(U\times V)\bigr).
\]
Since $(x,y)\in F(u,v)$, one has
\[
x\in F_{2}(v),\qquad y\in F_{1}(u).
\]
Moreover, using the product formula for the minimal time function and the
distance to the complement of the product set, the above inequality becomes
\[
\max\{\tau_{1},\tau_{2}\}\max\bigl\{T_{L}(x,u),\,T_{M}(y,v)\bigr\}<\min
\bigl\{\mathbf{d}(x,\mathbb{X}\setminus U),\,\mathbf{d}(y,\mathbb{Y}\setminus
V)\bigr\}.
\]
Consequently,
\[
\tau_{1}T_{L}(x,u)<\mathbf{d}(y,\mathbb{Y}\setminus V)
\]
and
\[
\tau_{2}T_{M}(y,v)<\mathbf{d}(x,\mathbb{X}\setminus U).
\]
Applying assumptions (i) and (ii), respectively, yields
\[
T_{M}\bigl(y,F_{1}(x)\bigr)\leq\tau_{1}T_{L}(x,u)
\]
and
\[
T_{L}\bigl(x,F_{2}(y)\bigr)\leq\tau_{2}T_{M}(y,v).
\]
Therefore,
\begin{align*}
T_{L\times M}\bigl((x,y),\mathcal{F}(x,y)\bigr) &  =\max\bigl\{T_{M}%
(y,F_{1}(x)),T_{L}(x,F_{2}(y))\bigr\}\\
&  \leq\max\bigl\{\tau_{1}T_{L}(x,u),\tau_{2}T_{M}(y,v)\bigr\}\\
&  \leq\max\{\tau_{1},\tau_{2}\}\max\bigl\{T_{L}(x,u),T_{M}(y,v)\bigr\}\\
&  =\max\{\tau_{1},\tau_{2}\}T_{L\times M}\bigl((x,y),(u,v)\bigr).
\end{align*}
Hence $F$ is directionally orbitally Aubin continuous on $U\times V$ w.r.t.
$L\times M$ with constant $\max\{\tau_{1},\tau_{2}\}$.  

\medskip

(b) We now consider the product space $\left(  \mathbb{X}\times\mathbb{Y}%
,\Vert\cdot\Vert_{3}\right)  $. Let $(x,y)\in U\times V$ and $(u,v)\in X\times
Y$ satisfy
\[
(x,y)\in\mathcal{F}(u,v)
\]
and
\[
\sqrt{\tau_{1}\tau_{2}}\,T_{\sqrt{\tau_{1}^{-1}}L\times\sqrt{\tau_{2}^{-1}}%
M}^{\Vert\cdot\Vert_{3}}\bigl((x,y),(u,v)\bigr)<\mathbf{d}%
\bigl((x,y),(\mathbb{X}\times\mathbb{Y})\setminus(U\times V)\bigr).
\]
Then
\[
x\in F_{2}(v),\qquad y\in F_{1}(u),
\]
and, by the definition of the norm $\Vert\cdot\Vert_{3}$,
\begin{align*}
&  \sqrt{\tau_{1}\tau_{2}}\max\Bigl\{\sqrt{\tau_{1}}\,T_{L}(x,u),\sqrt
{\tau_{2}}\,T_{M}(y,v)\Bigr\}\\
&  \qquad<\min\Bigl\{\sqrt{\tau_{1}}\,\mathbf{d}(x,\mathbb{X}\setminus
U),\sqrt{\tau_{2}}\,\mathbf{d}(y,\mathbb{Y}\setminus V)\Bigr\}.
\end{align*}
Hence,
\[
\tau_{1}T_{L}(x,u)<\mathbf{d}(y,\mathbb{Y}\setminus V)
\]
and
\[
\tau_{2}T_{M}(y,v)<\mathbf{d}(x,\mathbb{X}\setminus U).
\]
Invoking assumptions (i) and (ii), we obtain
\[
T_{M}(y,F_{1}(x))\leq\tau_{1}T_{L}(x,u)
\]
and
\[
T_{L}(x,F_{2}(y))\leq\tau_{2}T_{M}(y,v).
\]
Therefore,
\begin{align*}
&  T_{\sqrt{\tau_{1}^{-1}}L\times\sqrt{\tau_{2}^{-1}}M}^{\Vert\cdot\Vert_{3}%
}\bigl((x,y),\mathcal{F}(x,y)\bigr)\\
&  =\max\Bigl\{\sqrt{\tau_{1}}\,T_{L}(x,F_{2}(y)),\sqrt{\tau_{2}}%
\,T_{M}(y,F_{1}(x))\Bigr\}\\
&  \leq\max\Bigl\{\tau_{2}\sqrt{\tau_{1}}\,T_{M}(y,v),\tau_{1}\sqrt{\tau_{2}%
}\,T_{L}(x,u)\Bigr\}\\
&  \leq\sqrt{\tau_{1}\tau_{2}}\max\Bigl\{\sqrt{\tau_{1}}\,T_{L}(x,u),\sqrt
{\tau_{2}}\,T_{M}(y,v)\Bigr\}\\
&  =\sqrt{\tau_{1}\tau_{2}}\,T_{\sqrt{\tau_{1}^{-1}}L\times\sqrt{\tau_{2}%
^{-1}}M}^{\Vert\cdot\Vert_{3}}\bigl((x,y),(u,v)\bigr).
\end{align*}
This proves the claimed property and completes the proof. \hfill$\square$ 

\bigskip

%The next result shows that directional orbital Aubin continuity is preserved
%under the product construction introduced in the setting (S). The proof follows
%the same scheme as that of Proposition~\ref{pr2}, but requires a suitable
%adaptation to the orbital Aubin framework.

%We are now in a position to derive the announced coincidence and coupled fixed
%point results.
Now we are able to present the announced coincidence results. 

\begin{thm}
\label{TRON19} In the setting (S), let $U\subset X$ and $V\subset Y$ be open
sets, and let $\tau_{1},\tau_{2}>0$. Assume that: 

\begin{enumerate}
\item[(i)] $F_{1}$ is Milyutin directionally regular on $U\times V$ w.r.t. $L$
and $M$, as well as w.r.t. $-L$ and $-M$, with constant $\tau_{1}$; 

\item[(ii)] $F_{2}$ is Milyutin directionally regular on $V\times U$ w.r.t.
$M$ and $L$, as well as w.r.t. $-M$ and $-L$, with constant $\tau_{2}$; 

\item[(iii)] {\small
\[
\tau:=\max\{\tau_{1},\tau_{2}\}<1;
\]
}

\item[(iv)] there exists $(\overline{u},\overline{v})\in U\times V$ such that
one of the following conditions holds:
\begin{align*}
(\gamma)\qquad &  \max\Bigl\{T_{L}\bigl(\overline{u},F_{1}^{-1}(\overline
{v})\bigr),T_{M}\bigl(\overline{v},F_{2}^{-1}(\overline{u})\bigr)\Bigr\}\\
&  \hspace{2cm}<(1-\tau)\min\Bigl\{\mathbf{d}(\overline{u},\mathbb{X}\setminus
U),\mathbf{d}(\overline{v},\mathbb{Y}\setminus V)\Bigr\},\\[1ex]
(\delta)\qquad &  \max\Bigl\{T_{-L}\bigl(\overline{u},F_{1}^{-1}(\overline
{v})\bigr),T_{-M}\bigl(\overline{v},F_{2}^{-1}(\overline{u})\bigr)\Bigr\}\\
&  \hspace{2cm}<(1-\tau)\min\Bigl\{\mathbf{d}(\overline{u},\mathbb{X}\setminus
U),\mathbf{d}(\overline{v},\mathbb{Y}\setminus V)\Bigr\}.
\end{align*}

\end{enumerate}

Then the following assertions hold. 

\begin{enumerate}
\item[(a)] On the space $\left(  \mathbb{X}\times\mathbb{Y},\Vert\cdot
\Vert_{1}\right)  $, for every $\varepsilon>0$,
\[
\operatorname{Fix}(\mathcal{F})(\varepsilon)\cap(U\times V)\neq\emptyset.
\]
Moreover, if condition $(\gamma)$ holds, then
\[
\mathbf{d}\bigl((\overline{u},\overline{v}),\operatorname{Fix}(\mathcal{F}%
)(\varepsilon)\cap(U\times V)\bigr)\leq\frac{T_{L\times M}\bigl((\overline
{u},\overline{v}),\mathcal{F}^{-1}(\overline{u},\overline{v})\bigr)}{1-\tau},
\]
whereas if condition $(\delta)$ holds, then
\[
\mathbf{d}\bigl((\overline{u},\overline{v}),\operatorname{Fix}(\mathcal{F}%
)(\varepsilon)\cap(U\times V)\bigr)\leq\frac{T_{-(L\times M)}\bigl((\overline
{u},\overline{v}),\mathcal{F}^{-1}(\overline{u},\overline{v})\bigr)}{1-\tau}.
\]

\item[(b)] If $X$ and $Y$ are complete and both $F_{1}$ and $F_{2}$ have
closed graphs, then
\[
\operatorname{Fix}(\mathcal{F})\cap(U\times V)\neq\emptyset.
\]
Furthermore, if condition $(\gamma)$ holds, then
\[
\mathbf{d}\bigl((\overline{u},\overline{v}),\operatorname{Fix}(\mathcal{F}%
)\cap(U\times V)\bigr)\leq\frac{T_{L\times M}\bigl((\overline{u},\overline
{v}),\mathcal{F}^{-1}(\overline{u},\overline{v})\bigr)}{1-\tau},
\]
whereas if condition $(\delta)$ holds, then
\[
\mathbf{d}\bigl((\overline{u},\overline{v}),\operatorname{Fix}(\mathcal{F}%
)\cap(U\times V)\bigr)\leq\frac{T_{-(L\times M)}\bigl((\overline{u}%
,\overline{v}),\mathcal{F}^{-1}(\overline{u},\overline{v})\bigr)}{1-\tau}.
\]

\end{enumerate}
\end{thm}

\noindent Proof.  The conclusion follows directly from
Theorem~\ref{Main-Theorem} applied to the mapping $F$, together with
Proposition~\ref{pr2}(a). \hfill$\square$  

\begin{rmk}
The conclusions of part (b) ensure that for every $\eta>0$ there exists a pair
$(\hat{u},\hat{v})\in U\times V$ such that
\[
\hat{u}\in F_{2}(\hat{v}),\qquad\hat{v}\in F_{1}(\hat{u}),
\]
and
\[
\max\bigl\{\Vert\overline{u}-\hat{u}\Vert,\Vert\overline{v}-\hat{v}%
\Vert\bigr\}\leq(1+\eta)\frac{T_{L\times M}\bigl((\overline{u},\overline
{v}),\mathcal{F}^{-1}(\overline{u},\overline{v})\bigr)}{1-\tau},
\]
or alternatively,
\[
\max\bigl\{\Vert\overline{u}-\hat{u}\Vert,\Vert\overline{v}-\hat{v}%
\Vert\bigr\}\leq(1+\eta)\frac{T_{-(L\times M)}\bigl((\overline{u},\overline
{v}),\mathcal{F}^{-1}(\overline{u},\overline{v})\bigr)}{1-\tau}.
\]

\end{rmk}

\begin{rmk}
\begin{enumerate}

\item[(i)] By combining Theorem~\ref{Main-Theorem} with Proposition~\ref{pr2}%
(b), one obtains analogous coincidence results on the space $\left(
\mathbb{X}\times\mathbb{Y},\Vert\cdot\Vert_{2}\right)  $. Likewise,
Theorem~\ref{Direct Lip} together with Proposition~\ref{pr3} yields analogous
results on the product spaces $\left(  \mathbb{X}\times\mathbb{Y},\Vert
\cdot\Vert_{1}\right)  $ and $\left(  \mathbb{X}\times\mathbb{Y},\Vert
\cdot\Vert_{3}\right)  $. 

Moreover, global variants of Theorem~\ref{TRON19} are obtained by taking $U=X$
and $V=Y$. In this case, assumption (iv) is replaced by the requirement that
there exists $(\overline{u},\overline{v})\in X\times Y$ such that one of the
following conditions holds:
\begin{align*}
(\gamma)\qquad &  \max\Bigl\{T_{L}\bigl(\overline{u},F_{1}^{-1}(\overline
{v})\bigr),T_{M}\bigl(\overline{v},F_{2}^{-1}(\overline{u}%
)\bigr)\Bigr\}<\infty,\\[1ex]
(\delta)\qquad &  \max\Bigl\{T_{-L}\bigl(\overline{u},F_{1}^{-1}(\overline
{v})\bigr),T_{-M}\bigl(\overline{v},F_{2}^{-1}(\overline{u}%
)\bigr)\Bigr\}<\infty.
\end{align*}

\item[(ii)] In the framework of normed spaces, the above results recover, and
in several respects extend, recent coincidence theorems established in
\cite[Theorems~4.1--4.3]{T}. In comparison with related contributions such as
\cite{DF,I1,I3,ZO}, the assumptions employed here are substantially weaker,
while the resulting conclusions concerning coincidence points and coupled
fixed points remain of the same nature, differing essentially only in the
associated distance estimates. 
\end{enumerate}
\end{rmk}

\subsection{Applications to perturbation stability}

We now present an application of the preceding fixed point principles to the
stability of Milyutin-type regularity under directional perturbations. In
addition to the results obtained in \cite{I1}, \cite{I2}, and \cite{THN}, the
framework developed here also yields approximate solvability results in
possibly noncomplete spaces. 

\begin{thm}
\label{Fix-Milyu-Iof} Let $X$ and $Y$ be normed spaces, and let $U\subset X$
and $V\subset Y$ be open sets. Consider nonempty sets of directions
$\emptyset\neq L\subset S_{\mathbb{X}}$ and $\emptyset\neq M\subset
S_{\mathbb{Y}}$. Let $F:X\rightrightarrows Y$ be a set-valued mapping,
$g:X\rightarrow Y$ be a single-valued mapping, and let $\tau>0$. Assume that: 

\smallskip

\noindent(i) $F$ is Milyutin directionally regular on $U\times V$ w.r.t. $L$
and $M$ with constant $\tau$; 

\smallskip

\noindent(ii) $F$ is Milyutin directionally regular on $U\times V$ w.r.t. $-L$
and $-M$ with the same constant $\tau$; 

\smallskip

\noindent(iii) there exists $\lambda\in(0,\tau^{-1})$ such that, for every
$u\in U$ and every $v\in u+\operatorname{cone}L$,
\[
T_{-M}(g(u),g(v))\leq\lambda T_{L}(u,v),
\]
and, for every $u\in U$ and every $v\in u-\operatorname{cone}L$,
\[
T_{M}(g(u),g(v))\leq\lambda T_{-L}(u,v);
\]

\smallskip

\noindent(iv) there exists an open set $W\subset Y$ such that
\[
W-g(U)\subset V.
\]

Let $(x,y)\in U\times W$ satisfy
\[
0<(\tau^{-1}-\lambda)^{-1}T_{M}(y,F(x)+g(x))<\mathbf{d}(x,\mathbb{X}\setminus
U),
\]
and define the set-valued mapping
\[
R_{y}:=F^{-1}\circ(-g+y).
\]

Then the following assertions hold. 

\smallskip

\noindent(a) For every $\varepsilon>0$,
\[
\operatorname{Fix}(R_{y})(\varepsilon)\cap U\neq\emptyset,
\]
and
\[
\mathbf{d}\bigl(x,\operatorname{Fix}(R_{y})(\varepsilon)\cap U\bigr)\leq
\frac{\tau}{1-\tau\lambda}\,T_{M}(y,F(x)+g(x)).
\]

\smallskip

\noindent(b) If $X$ is Banach, $\operatorname{Gr}F$ is closed and $g$ is
continuous, then
\[
\operatorname{Fix}(R_{y})\cap U\neq\emptyset,
\]
and
\[
\mathbf{d}\bigl(x,\operatorname{Fix}(R_{y})\cap U\bigr)\leq\frac{\tau}%
{1-\tau\lambda}\,T_{M}(y,F(x)+g(x)).
\]
Equivalently,
\[
\mathbf{d}\bigl(x,(F+g)^{-1}(y)\cap U\bigr)\leq\frac{\tau}{1-\tau\lambda
}\,T_{M}(y,F(x)+g(x)).
\]

\end{thm}

\noindent Proof.  Fix $(x,y)$ as in the assumptions. We verify that the
mapping $R_{y}$ satisfies the hypotheses of Theorem~\ref{Direct Lip}. 

We first prove that $R_{y}$ is directionally orbitally Aubin continuous on $U$
w.r.t. $L$ with constant $\tau\lambda$. Let $x^{\prime}\in U$ and $u^{\prime
}\in X$ satisfy
\[
x^{\prime}\in R_{y}(u^{\prime})
\]
and
\[
\tau\lambda T_{L}(x^{\prime},u^{\prime})<\mathbf{d}(x^{\prime},\mathbb{X}%
\setminus U).
\]
Since $x^{\prime}\in R_{y}(u^{\prime})$, one has
\[
y-g(u^{\prime})\in F(x^{\prime}).
\]
Moreover, the condition $T_{L}(x^{\prime},u^{\prime})<+\infty$ implies that
\[
u^{\prime}\in x^{\prime}+\operatorname{cone}L.
\]
By assumption~(iv),
\[
y-g(x^{\prime})\in W-g(U)\subset V.
\]
Using the first inequality in assumption~(iii), we obtain
\begin{align*}
\tau T_{M}(y-g(x^{\prime}),F(x^{\prime})) &  \leq\tau T_{M}(y-g(x^{\prime
}),y-g(u^{\prime}))\\
&  =\tau T_{-M}(g(x^{\prime}),g(u^{\prime}))\\
&  \leq\tau\lambda T_{L}(x^{\prime},u^{\prime})\\
&  <\mathbf{d}(x^{\prime},\mathbb{X}\setminus U).
\end{align*}
Since $F$ is Milyutin directionally regular on $U\times V$ w.r.t. $L$ and $M$
with constant $\tau$, it follows that
\begin{align*}
T_{L}(x^{\prime},R_{y}(x^{\prime})) &  =T_{L}\bigl(x^{\prime-1}(y-g(x^{\prime
}))\bigr)\\
&  \leq\tau T_{M}(y-g(x^{\prime}),F(x^{\prime}))\\
&  \leq\tau\lambda T_{L}(x^{\prime},u^{\prime}).
\end{align*}
Hence $R_{y}$ is directionally orbitally Aubin continuous on $U$ w.r.t. $L$
with constant $\tau\lambda$. 

We next show that $R_{y}$ is directionally orbitally Aubin continuous on $U$
w.r.t. $-L$ with the same constant $\tau\lambda$. Let $x^{\prime}\in U$ and
$u^{\prime}\in X$ satisfy
\[
x^{\prime}\in R_{y}(u^{\prime})
\]
and
\[
\tau\lambda T_{-L}(x^{\prime},u^{\prime})<\mathbf{d}(x^{\prime},\mathbb{X}%
\setminus U).
\]
Since
\[
u^{\prime}\in x^{\prime}-\operatorname{cone}L,
\]
the second inequality in assumption~(iii) yields
\begin{align*}
\tau T_{-M}(y-g(x^{\prime}),F(x^{\prime})) &  \leq\tau T_{-M}(y-g(x^{\prime
}),y-g(u^{\prime}))\\
&  =\tau T_{M}(g(x^{\prime}),g(u^{\prime}))\\
&  \leq\tau\lambda T_{-L}(x^{\prime},u^{\prime})\\
&  <\mathbf{d}(x^{\prime},\mathbb{X}\setminus U).
\end{align*}
Using the Milyutin directional regularity of $F$ on $U\times V$ w.r.t. $-L$
and $-M$, we obtain
\begin{align*}
T_{-L}(x^{\prime},R_{y}(x^{\prime})) &  =T_{-L}\bigl(x^{\prime-1}%
(y-g(x^{\prime}))\bigr)\\
&  \leq\tau T_{-M}(y-g(x^{\prime}),F(x^{\prime}))\\
&  \leq\tau\lambda T_{-L}(x^{\prime},u^{\prime}).
\end{align*}

Consequently, assumption~(i) of Theorem~\ref{Direct Lip} is satisfied for the
mapping $R_{y}$. 

We now verify the second hypothesis of Theorem~\ref{Direct Lip}. Since
\[
T_{M}(y,F(x)+g(x))<(\tau^{-1}-\lambda)\mathbf{d}(x,\mathbb{X}\setminus U),
\]
we have
\[
\tau T_{M}(y,F(x)+g(x))<(1-\tau\lambda)\mathbf{d}(x,\mathbb{X}\setminus U).
\]
Observe also that $y-g(x)\in V$. Hence, by the Milyutin directional regularity
of $F$ on $U\times V$ w.r.t. $L$ and $M$,
\begin{align*}
T_{L}(x,R_{y}(x)) &  =T_{L}\bigl(x,F^{-1}(y-g(x))\bigr)\\
&  \leq\tau T_{M}(y-g(x),F(x))\\
&  =\tau T_{M}(y,F(x)+g(x))\\
&  <(1-\tau\lambda)\mathbf{d}(x,\mathbb{X}\setminus U).
\end{align*}
Therefore all assumptions of Theorem~\ref{Direct Lip} are fulfilled. 

\smallskip

\noindent(a) Applying Theorem~\ref{Direct Lip}, we conclude that, for every
$\varepsilon>0$,
\[
\operatorname{Fix}(R_{y})(\varepsilon)\cap U\neq\emptyset,
\]
and
\[
\mathbf{d}\bigl(x,\operatorname{Fix}(R_{y})(\varepsilon)\cap U\bigr)\leq
\frac{T_{L}(x,R_{y}(x))}{1-\tau\lambda}.
\]
Using the preceding estimate for $T_{L}(x,R_{y}(x))$, we deduce that
\[
\mathbf{d}\bigl(x,\operatorname{Fix}(R_{y})(\varepsilon)\cap U\bigr)\leq
\frac{\tau}{1-\tau\lambda}\,T_{M}(y,F(x)+g(x)).
\]

\smallskip

\noindent(b) Assume now that $X$ is Banach, $\operatorname{Gr}F$ is closed and
$g$ is continuous. Then $\operatorname*{Gr}R_{y}$ is closed and another
application of Theorem~\ref{Direct Lip} yields
\[
\operatorname{Fix}(R_{y})\cap U\neq\emptyset,
\]
together with
\[
\mathbf{d}\bigl(x,\operatorname{Fix}(R_{y})\cap U\bigr)\leq\frac{T_{L}%
(x,R_{y}(x))}{1-\tau\lambda}\leq\frac{\tau}{1-\tau\lambda}\,T_{M}%
(y,F(x)+g(x)).
\]

Finally, observe that
\[
u\in\operatorname{Fix}(R_{y})\iff u\in R_{y}(u)\iff y\in F(u)+g(u)\iff
u\in(F+g)^{-1}(y).
\]
Hence
\[
\operatorname{Fix}(R_{y})=(F+g)^{-1}(y),
\]
and therefore
\[
\mathbf{d}\bigl(x,(F+g)^{-1}(y)\cap U\bigr)\leq\frac{\tau}{1-\tau\lambda
}\,T_{M}(y,F(x)+g(x)).
\]
The proof is complete. \hfill$\square$ 

\begin{rmk}
The assumptions (ii)--(iii) in Theorem~\ref{Fix-Milyu-Iof} may be replaced by
the following variants: 

\smallskip

\noindent(ii)$^{\prime}$ $F$ is Milyutin directionally regular on $U\times V$
w.r.t. $-L$ and $M$ with constant $\tau$; 

\smallskip

\noindent(iii)$^{\prime}$ there exists $\lambda\in(0,\tau^{-1})$ such that,
for every $u\in U$ and every $v\in u+\operatorname{cone}L$,
\[
T_{-M}(g(u),g(v))\leq\lambda T_{L}(u,v),
\]
and, for every $u\in U$ and every $v\in u-\operatorname{cone}L$,
\[
T_{-M}(g(u),g(v))\leq\lambda T_{-L}(u,v).
\]

\end{rmk}

\begin{rmk}
\begin{enumerate}

\item Local counterparts of Theorem~\ref{Fix-Milyu-Iof} can be obtained under
additional assumptions, for instance when $\operatorname{cone}L$ and
$\operatorname{cone}M$ are convex. 

\item Global versions of Theorem~\ref{Fix-Milyu-Iof} follow by taking $U=X$
and $V=Y$ and assuming that $(x,y)\in X\times Y$ satisfies
\[
0<T_{M}(y,F(x)+g(x))<+\infty.
\]

\end{enumerate}
\end{rmk}

\begin{rmk}
Theorem~5.2 in \cite{T} is recovered from Theorem~\ref{Fix-Milyu-Iof} by
taking $L=S_{\mathbb{X}}$ and $M=S_{\mathbb{Y}}$. 
\end{rmk}

\section{Concluding remarks}

In this paper, we developed a directional framework for the study of fixed
point phenomena of set-valued mappings based on directional orbital regularity
and directional orbital Aubin continuity. The obtained results extend several
recent fixed point and stability theorems to a genuinely directional setting
formulated through minimal time functions associated with prescribed sets of
admissible directions. 

Within this framework, we established approximate and exact fixed point
theorems under weak directional assumptions and derived corresponding
coincidence point and stability results. The proposed approach reveals a
nontrivial interplay between directional regularity properties and orbital
iterative constructions, allowing one to capture asymmetric and anisotropic
behaviors that are not accessible through the classical nondirectional theory. 

The results obtained herein suggest several directions for future research. In
particular, it would be of interest to further investigate the influence of
the admissible direction sets and their geometric properties on the structure
of fixed point and coincidence point sets, as well as possible applications to
vector optimization, generalized equations, and variational systems governed
by set-valued mappings. 

%In this work, we have refined and extended several recent fixed point results
%to a directional setting. Furthermore, coincidence and stability results are
%also obtained within a fully directional framework, bringing novel insights
%and a deeper understanding of some phenomena that cannot be observed in the
%classical setting. The applications we have provided offer new perspectives
%for further study. In saying this, we have in mind a deeper understanding of
%the role of the initially considered directions and their opposites as they
%appear in the main results of this paper, with possible applications to vector
%optimization problems governed by set-valued maps.

\bigskip

\noindent Data availability. This manuscript has no associated data. 

\noindent Disclosure statement. The authors contributed equally to all aspects
of the study and to the preparation of the manuscript. No potential conflict
of interest was reported by the authors. 

\noindent Funding. The authors declare that no funds, grants, or other support
were received during the preparation of this manuscript 

\bigskip
\bibliographystyle{plain}
\bibliography{TDT.bib}

\end{document}